\documentclass[11pt]{article}

\usepackage{amsmath,amssymb,vmargin}
\usepackage[applemac]{inputenc}

\newcommand{\R}{\mathbb R}

\newcommand{\Z}{\mathbb Z}

\newcommand{\eps}{\varepsilon}

\newcommand{\sumZe}{\sum_{n \ne 0}}
\newcommand{\ud}{u_D}

\renewcommand{\d}{\partial}
\newcommand{\DtoN}[1]{\Lambda_{#1}}

\newcommand{\dom}{\Omega}

\newcommand{\bdom}{\d\dom}

\newcommand{\sH}{{\rm H}}
\newcommand{\sW}{{\rm W}}

\renewcommand{\div}[1]{\mbox{div}\left(#1\right)}

\newtheorem{thm}{Theorem}[section]
\newtheorem{prop}[thm]{Proposition}
\newtheorem{lemme}[thm]{Lemma}

\newtheorem{remark}[thm]{Remark}

\def\blacksquare{
\thinspace\nobreak \vrule width 5pt height 5pt depth 0pt}

\newenvironment{proofof}[1]{\begin{trivlist}
\item[]\hspace{0cm}{\bf Proof of #1: }
\hspace{0cm} }{\hfill $\blacksquare$
\end{trivlist}}

\begin{document}

\title{A remark on precomposition on $\sH^{1/2}(S^1)$ and $\eps$-identifiability of disks in tomography.}
\author{M. Dambrine and D. Kateb \footnote{LMAC Universit\'e de Technologie de Compi\`egne BP 20529 60205 Compi\`egne Cedex -  France. ~~~~~~~~~ Email : Marc.Dambrine@utc.fr, Djalil.Kateb@utc.fr}}

\maketitle

\begin{abstract}
We consider the inverse conductivity problem with one measurement
for the equation $div((\sigma_1+(\sigma_2-\sigma_1)\chi_D)\nabla{u})=0$ determining the unknown inclusion $D$ included in $\Omega$. We suppose that
$\Omega$ is the unit disk of $\mathbb{R}^2$. With the tools of the
conformal mappings, of elementary Fourier analysis and also the
action of some quasi-conformal mapping on the Sobolev space
$\sH^{1/2}(S^1)$, we show how to approximate the
Dirichlet-to-Neumann map when the original inclusion $D$ is a
$\varepsilon-$ approximation of a disk. This enables us to give
some uniqueness and stability results.
\end{abstract}



\section{Introduction.}

In this paper, we study the inverse problem of conductivity with
one measurement. Given a bounded domain $\dom\subset \mathbb{R}^2$
with reasonably smooth boundary, an connected open set $D$ compactly
contained in $\Omega$, we consider for any $f \in \sH^{1/2}(\partial D)$ the problem of recovering the subset $D$
entering the Dirichlet equation
\begin{equation}
\label{equation:conductivite} P[D,f]\left\{
\begin{array}{rlll}
\div{\left(\sigma_{1}+ (\sigma_{2}-\sigma_{1}) \chi_{D} \right)
\nabla u} &=& 0,~ &\text{ in } \dom,\\
u&=&f,~&\text{ on } \partial \dom
\end{array}
\right.
\end{equation}
from the knowledge of the current flux $g=\sigma_{1} \d_n u$ in
$\bdom$ induced by the boundary value $f=u_{|\bdom}$. We will
denote $\Lambda_{D} : \sH^{1/2}(\bdom)\rightarrow
\sH^{-1/2}(\bdom)$ the Dirichlet-to-Neumann map which maps the
Dirichlet data $f$ onto the corresponding Neumann data
$g=\Lambda_{D}(f)=\sigma_{1} \d_n u $. We can reformulate the
inverse problem with one measurement as the determination of $D$
from the Cauchy pair $(f,g)$. We mention here that we do not need
the full knowledge of the Dirichlet-to-Neumann map but only one
pair of Cauchy data $(f,g)$. For such a problem, we know that the
uniqueness question is, in general, an open problem. It has been
solved only for the special class of convex polyhedra, disks and
balls. For other domains, Fabes, Kang and Seo \cite{FabKanSeo} have studied the
global uniqueness and stability within the class of domains which
are $\eps-$ perturbations of disks. The main ingredients in the work were layer potential techniques
and representation formula for the solution $u_D$ of the problem $P[D,f]$.

Our main goal is to revisit the paper \cite{FabKanSeo} with other techniques than boundary integral
representations. Throughout our paper, the two-dimensional case will be considered. Instead of layer potential techniques,  conformal mappings and Fourier analysis will be another approach
to review the two questions of stability and uniqueness within the class of disks and
perturbed disks.

Let us illustrate briefly the main steps of our arguments.
Since $\Omega$ is doubly connected,  conformal mappings allows the construction of  the conformal transplant function $u$ which is solution of an elliptic problem that is obtained by transporting the original problem $P[D,f]$ by means of a change of variables induced by the conformal mappings.  A natural way to study the original Dirichlet-to-Neumann map is to study the transplanted Dirichlet-to-Neumann map; indeed one can show that when the original $D$ is a disk in $\Omega$, then we can give an expression of the new Dirichlet-to-Neumann map by means of Moebius transforms and the classical formula of the Dirichlet-to -Neumann operator related to an concentric annulus. The elementary properties of Moebius transforms allow us to get an  uniqueness result within the class of circular inclusions and for some special Dirichlet boundary measurement.

When $D$ is not a disk, things become more difficult. The conformal transplantation furnishes a Dirichlet-to-Neumann operator that is not very convenient to study. Indeed,  all the classical tools of perturbation theory  by compact operators and  Von-Neumann expansion of the inverses are absent  and thus  we have no way to make our explicit formula more suitable for numerical purposes. However, when the original inclusion $D$ is an $\varepsilon-$ perturbation of a disk $B$, then one can show that we have a reliable expression of the form  $\Lambda_D=\Lambda_B+R_\varepsilon$
with a remainder $R_\varepsilon$ that is of order $\varepsilon^\alpha$. We show that   $\alpha$ depends on the Sobolev regularity of the  Dirichlet boundary measurement $f$ and of the regularity of the boundary $\partial D$.
In the conformal  transplant, we have to deal with the two conformal mappings that map respectively $\Omega \backslash \overline{D}$ into the annulus and and $D$ on the ball;  and  the restriction of the maps on the corresponding boundaries will be of great importance in the error estimate. In our context, the diffeomorphism $\xi$  is obtained from a composition of the two  boundary correspondence functions. The error estimate is not straighforward, it  is a consequence of the hardest  problem   of estimating $\| h - h \circ \xi\|_{\sH^{1/2}(S^1)}$ when $\xi : S^1 \mapsto S^1 $
is a $\sW^{1,\infty}$ diffeomorphism of the circle and when $h$ is a function that belongs to some Sobolev space $\sH^s(S^1)$,  $s>{1\over 2}$. We were not able to give the best  Sobolev exponent $s$ for which the estimate is true. However we give a result for the exponent values  $s=1+\alpha$ for some $0<\alpha<1.$ At our best knowledge, the question remains open when $h$ belongs to  $\sH^s(S^1)$ when ${1\over 2}<s<1$.
Our  result about the precomposition of Sobolev spaces with quasi-regular diffeomorphisms are in the  continuation of the pioneering works (see \cite{BouBreMir,Bourdaud}) where are studied   the action of  quasi-regular homeomorphisms on the critical Sobolev space $\sH^{1/2}({S^1})$.

Let us point out that the resolution of an inverse boundary value problem for harmonic functions arising in electrostatic imaging through conformal mapping techniques has been introduced by Kress and his
collaborators. The interested  reader  can consult the seminal work of Kress and al \cite{AkdKre,HadKre} and our paper in \cite{DamKat}.

The paper is organized as follows. In section 2, after
introducing some definitions and recalling some preliminary
results concerning Moebius conformal mapping, we state the
uniqueness results for disks. In section 3, we investigate the
continuity properties of the superposition operators on
$\sH^{1/2}(S^1)$ generated by  regular diffeomorphisms of the
circle. We then describe the approximation of the
Dirichlet-to-Neumann map obtained after a sufficiently small
deformation of a disk. In section 4, we prove the main result of
uniqueness for disk and the $\eps$ identifiability of $\eps$ disks.
In section 5, we prove the  precomposition inequality.

%%%%%%%%%%%%%%%%%%%%%%%%%%%%%%
\section{Main assumptions and results}
%%%%%%%%%%%%%%%%%%%%%%%%%%%%%%

We shall assume throughout that $\Omega$ is the unit ball of
$\mathbb{R}^2$. Let us introduce the notion of small perturbation of
disks. Given $\varepsilon \ge 0$, a $C^2$ domain $D$ is called an
$\varepsilon-$ perturbation of a disk if there exists $\delta \in
C^2(\partial B)$ with $\| \delta \|_{C^2(\partial B)}
<1$ such that
\begin{equation*}
\partial D : x+\varepsilon \delta(x)\nu(x),~~~x\in\partial B
\end{equation*}
where $\nu(x)$ is the outward unit normal to $\partial B$ at $x$.
Denoting $\Omega_0 \subset \dom$ the set of points at some distance $\delta_0$ from
$\partial \dom$, we will denote by $C[\varepsilon]$ the class of
$\varepsilon-$ perturbations of all disks contained in $\Omega_0$
with the radius larger than a fixed number $\rho_0$ that can be
arbitrary small provided than it remains big with respect to $\eps$.
\par
\noindent We will assume that the domain $D$ entering in equation
\eqref{equation:conductivite}  is a disk or an $\varepsilon-$ perturbation of a disk
$B\subset \Omega_0$. Our main results concern the
$identifiability$ (the case of a perfect disk) and the approximate
identifiability.

In a first time, we deal with the perfect case where $\varepsilon=0$; we
have
\begin{thm}
\label{unicite:disques} Let $D_1$ be a disk centered at the origin
and of radius $R_1$. Then there exists a boundary Dirichlet
measurement $f(\theta)=\cos{\theta}$ such that if $D_2$ is an
arbitrary disk contained in $\dom$ and $\Lambda_{D_1}(f)
=\Lambda_{D_2}(f)$ then $D_1=D_2$.
\end{thm}
In a second time, we consider the case of perturbed disks. We take
the same boundary measurment $f(\theta)=\cos{\theta}$. We then
have the following
\begin{thm}
\label{unicite:fabes}
Let $D_0 \in C[\varepsilon]$ and let
$\Lambda_{D_0}(f)=g$. Then there exists a positive constant $C>0$
such that if $D \in C[\varepsilon]$ and $\Lambda_D(f)=g$ on
$\partial \dom$, then
\begin{equation}
| D\Delta D_0 | \le C \varepsilon^\alpha,
\end{equation}
where $0<\alpha<1$ is a constant depending only on the a priori
data and where $D\Delta D_0$ denotes the symmetric difference of
$D$ and $D_0$.
\end{thm}
\par
\noindent The Sobolev $\sH^{1/2}(S^1)$ plays an
important role in the proof of stability. We recall that $\sH^{1/2}(S^1)$ 
stands for the Hilbert space of real functions $f$ defined
on $S^1$ (modulo the constants) whose Fourier expansion
$$f(e^{i\theta})=\sum_{n = -\infty}^{n=\infty}c_n(f) e^{in\theta}$$
where the Fourier coefficients $(c_n(f))$ are such that the
sequence $(\sqrt{n}c_n(f))_n$ is square summable. For each $f$
belonging to the Sobolev space $\sH^{1/2}(S^1)$, its norm is the weighted $l^2$
norm $\left( \sum_{n=-\infty}^{\infty} | n| | c_n(f)|^2
\right)^{{1/ 2}}$. We will also write that $f$ belongs to the
half Sobolev space if and only if the sequence of its Fourier
transform $(c_n(f))$ belongs to $l_2^{{1 / 2}}(\mathbb{Z})$.

Let us recall some results about the action or composition (it is
en fact a ''precomposition'') by quasi-symmetric homeomorphisms of
the circle $S^1$. Given an orientation preserving homeomorphism
$\xi : S^1\mapsto S^1$ of the circle, we consider the superposition
operator $F_\xi$ generated by $\xi$ defined by
\begin{equation*}
F_\xi(f)=f\circ\xi,~f\in \sH^{{1/2}}(S^1).
\end{equation*}
It is known that $F_\xi$ maps $\sH^{1 /2}(S^1)$ onto itself if and
only if $\xi$ is quasi-symmetric in the sense that we must have
the doubling condition
\begin{equation*}
\cfrac{| \xi(2I)|}{| \xi(I)|}\le K,
\end{equation*}
where $K>0$ is positive, where $I$ is any interval on $S^1$ of
length less than $\pi$ and where $2I$ is the interval of $S^1$
after doubling $I$ but by keeping the same midpoint. Furthermore,
we have
\begin{equation*}
\| F_\xi\|_{{\cal L}(\sH^{{1/ 2}}(S^1),\sH^{{1/ 2}}(S^1))}\le \sqrt{K+{1\over K}}.
\end{equation*}
We recall also  that among all quasi-symmetric homeomorphisms of
$S^1$, the Moebius transformations of $S^1$ act unitarily on
$\sH^{1/ 2}(S^1)$.

An important question arises  : can we hope to bound the error
norm $\| F_\xi(u) -u\|_{\sH^{1/ 2}(S^1)}$ by $C \| u \|_{\sH^{1/
2}(S^1)}$.  The answer is important  since it will allow  us to
estimate the error between the original Dirichlet-to-Neumann
operator $\Lambda_D$ and the \emph{transformed} $\Lambda_B$. A
first idea is to guess  that such an estimate   can be  possible
if $\xi$ is not \emph{far} from a Moebius transform of $S^1$.
However, at this stage of our work, we have to add  some
regularity assumptions on the target function $u$. To be more
precise, we are only able to prove the following result of continuity 
for the precomposition by diffeomorphisms.

\begin{thm}
\label{theoreme:precomposition} Let $0 < \delta <1$ and let   $u$
be a function  belonging  to $\sH^{1+\delta}$. Let $\phi$ be a
regular  quasi-regular function that we suppose  to be a
$\sW^{1,\infty}$ diffeomorphism on $S^1$. Then we have
\begin{equation}\label{inegalite:composition}
\|u\circ \phi - u\|_{\sH^{1/2}(S^1)} \leq C(\delta')~
\|u\|_{\sH^{1+\delta}(S^1)} ~\omega_{\delta'}(\|\phi-I\|_{\sW^{1,\infty}(S^1)})
\end{equation}
holds for all $\delta' \in (0,\delta)$ where $\omega_{\delta'}$ is
the modulus of continuity defined by
\begin{equation}\label{definition:omega}
\omega_{\delta'} (t) =\max(t^{\delta'+1/2},t^{\delta'} ).
\end{equation}
\end{thm}

%%%%%%%%%%%%%%%%%%%%%%%%%%%%%%%%%%%%%
\section{Change of variables and superposition operators.}
\label{section:change_of_variables_and_superposition_operators}
%%%%%%%%%%%%%%%%%%%%%%%%%%%%%%%%%%%%%

We suppose that $D \in C[\varepsilon]$ is an $\varepsilon-$
perturbation of a disk. We aim to  approximate  the
Dirichlet-to-Neumann operator $\Lambda_D$ by $\Lambda_B$. In a
first time, we wish to transport the original problem $P[D,f]$ by
means of conformal transforms and to study the corresponding
Dirichlet-to-Neumann operator.
\subsection{Change of variables
and analysis of the Dirichlet-to-Neumann map} Thanks to the
classical mapping theorems, we know that there exists
$\rho\in(0,1)$ and an analytic function $\Phi^e$ that maps
bijectively $\dom\setminus\overline{D}$ onto the annulus
$\dom\setminus\overline{B}_\rho$ where $B_\rho$ is the centered
disk of radius $\rho$. If the outer boundaries correspond to each
other and if the image of one point on $\d\dom$ is prescribed then
$\Phi^e$ is uniquely determined. Furthermore, thanks to the
Riemann's mapping theorem, we know that there exists also  a
conformal mapping $\Phi^i$ that maps bijectively $D$ onto
$B_\rho$. We recall that the restrictions of the
conformal maps to the inclusion have the same regularity than
$\partial D$ (see \cite{Henrici1,Henrici3} for more
details).

We denote by $\Psi^i=(\Phi^i)^{-1}$ (respectively $(\Psi^e)^{-1})$
the inverse of $\Phi^i$ (respectively $\Phi^e$) and by $\gamma:
[0,|\d D|] \rightarrow \d D$ the parametrization of $\d D$ in
terms of arc-length. We set
\begin{equation}\label{definition:phi^i}
\phi^i(\theta) = \gamma^{-1}\left(\Psi^i(\rho e^{i\theta} )
\right),
\end{equation}
and
\begin{equation}\label{definition:phi^e}
\phi^e(\theta) = \gamma^{-1}\left(\Psi^e(\rho e^{i\theta} )
\right).
\end{equation}
Let $U_{B_\rho}^e$ (respectively $U_{B_\rho}^i$) denote the
conformal transplant of $(\ud)_{|\dom\setminus\overline{D}}$
(respectively $(\ud)_{|D}$). We have
\begin{equation*}
u^e=U_{B_\rho}^e\circ\Phi^e
\end{equation*}
and
\begin{equation*}
u^i=U_{B_\rho}^i\circ\Phi^i;
\end{equation*}
let us  give the explicit form of the elliptic equations
satisfied  each of the conformal transplants. We have

\begin{prop}
\label{proposition:transport}
Let $\xi$ be the diffeomorphism on $\d B_\rho$ defined by $\xi =
(\phi^e)^{-1}\circ \phi^i$. Then, we have
\begin{eqnarray*}
\Delta U_{B_\rho}^e &=& 0 \text{ in } \dom\setminus\overline{B}_\rho, \\
U_{B_\rho}^e &=& f \circ\Psi^e \text{ on } \d\dom,\\
\Delta U_{B_\rho}^i &=& 0 \text{ in } B_\rho,\\
U_{B_\rho}^i &=& U_{B_\rho}^e\circ \xi \text{ on } \d B_\rho,\\
\sigma_1 \left( \d_r U_{B_\rho}^e\circ \xi \right) \xi' &=&
\sigma_2 \d_r U_{B_\rho}^i \text{ on }\d B_\rho.
\end{eqnarray*}
\end{prop}

\begin{proofof}{Proposition \ref{proposition:transport}} It is elementary and essentially based on the
Cauchy-Riemann identities. We left the details to the reader.
\end{proofof}

\subsection{Boundary correspondance for $\varepsilon$ perturbations of disks.}

The important fact to notice now is that if $D \in
C[\varepsilon]$, then $\xi$ is a perturbation of the identity.
More precisely, one has the following result
\begin{prop}\label{proposition:petite:taille:xi}
There exists a positive constant $C>0$ such that 
\begin{equation}
\label{petite:taille:xi} \| \xi-I \|_{W^1_{\infty}(S^1)}\le C
\varepsilon
\end{equation}
holds for  all $D\in C[\varepsilon]$.
\end{prop}

\begin{proofof}{Proposition \ref{proposition:petite:taille:xi}}
Since $\xi = (\phi^e)^{-1}\circ \phi^i$, the proof is split into
two parts: in a first time we estimate the contribution of the
interior then in a second time the contribution of the exterior.
We claim that
\begin{equation}
\label{claim1} \| \phi^i-I \|_{W^1_{\infty}(S^1)}\le C \varepsilon,
\end{equation}
and
\begin{equation}
\label{claim2} \| \phi^e-I \|_{W^1_{\infty}(S^1)}\le C \varepsilon.
\end{equation}
Deducing \eqref{petite:taille:xi} from the claims \eqref{claim1},\eqref{claim2} is easy and left to the reader. We now prove the claims in the two next sections.
\end{proofof}

\subsubsection{Proof of claims \eqref{claim1}-\eqref{claim2}.}
\paragraph{The simply  connected case: claim \eqref{claim1}.}
Without loss of generality, one can  assume $\partial
\omega$ to be starlike with respect  to the origin. We use polar
coordinates to write
\begin{equation}
\partial D : ~~z=z(\phi)=r(\phi)e^{i\phi},~0 \le \phi \le 2\pi,
\end{equation}
where $r$ is a given positive regular function of period $2\pi$
such that
\begin{equation}\label{rayonpert}
r(\theta)=(R-\delta(\theta)),~0\le \theta \le 2\pi
\end{equation}
 $\delta$ being  a function of period $2\pi$ satisfying
$$
\sup_{0\le \theta<2\pi}| \delta^{(k)}|
<\varepsilon,~k=0,1,2.
$$
From Henrici (\cite{Henrici1},\cite{Henrici3}), we learn that  $\theta$ and $\phi^i(\theta)$ are
related by the Theodersen's integral equation
\begin{equation}\label{theodersen}
\phi^i(\theta)-\theta={\cal H}(\log{r(\phi(\theta))})
\end{equation}
where $\cal{H}$, the Hilbert transform on the circle, is  defined by
\begin{equation}
{\cal H}f(\theta)={1\over 2\pi}P.V.
\int_0^{2\pi}f(t)\cot{{\theta-t\over 2}}~dt.
\end{equation}
The same author learns us that the Theordersen's integral equation
admits exactly one continuous solution $\phi(\theta),~0\le \theta
\le 2\pi$ under the condition that the ratio
\begin{equation}\label{condition} \delta=\sup_{0\le\phi\le
2\pi}\left| {r'(\theta)\over r(\theta)}\right|
\end{equation}
satisfies  $\delta <1$. This  condition \eqref{condition} means
that the angle between the outward normal and the radius vectors
does not exceed $\arctan{\delta}<\delta$. It also means that we
have to deal with curves $\partial \omega$ that are not too far
from a circle. It is referred  as the $\delta$ condition. In our
context, $D$ is not to far from a disk.  Some straightforward
arguments (primarily due to Montel and Lindelof) show that the
boundary correspondence between $\theta$ and $\phi^i(\theta)$ is
given by
\begin{equation}
\phi^i(\theta)=\theta-{\cal H}\delta(\theta)+O(\varepsilon^2).
\end{equation}
The interested reader will  find a geometric proof in
(\cite{Henrici1}). Hence, if the perturbation $\delta$  is in
$W^2_{\infty}$, one can show easily that there exists a positive
constant $C>0$ such that
\begin{equation}
| \phi^i(\theta) -\theta| < C \varepsilon.
\end{equation}

\paragraph{The doubly connected case: claim \eqref{claim2}.}
As we did in the previous paragraph, we give the asymptotic
behavior of $\phi^e-I$  when $\varepsilon \rightarrow 0$. The
analog of the Theodersen's equations for  the doubly connected
case is described by the so called Theodersen's and Garrick
equations. The boundary correspondence is given by the following
result.

\begin{thm}
Let $O$ be a doubly connected region conformally  equivalent to
the annulus $\rho <| w| <1$. We suppose $O$ bounded by
piecewise analytic curves $\Gamma_0$ and $\Gamma_1$, both starlike
with respect to the origin and parametrized as above. If we
suppose that
\begin{equation}
\int_0^{2\pi}(\phi_0(\theta)-\theta)~d\theta=\int_0^{2\pi}(\phi_1(\theta)-\theta)~d\theta=0
\end{equation}
then there holds
\begin{equation}
\left\{
\begin{array}{lll}
\phi_0(\theta)-\theta&=&-{\cal K}_\rho(\log{r_1(\phi_1(\theta))})\\
\phi_1(\theta)-\theta&=&-{\cal H}_\rho(\log{r_1(\phi_1(\theta))})
\end{array}
\right.
\end{equation}
where the operators  ${\cal K}_\rho : L^2((0,2\pi))\mapsto
L^2((0,2\pi))$ and ${\cal H}_\rho : L^2((0,2\pi))\mapsto
L^2((0,2\pi)) $ are defined as follows
\begin{equation}
e^{im\theta}\mapsto {\cal H}_\rho(e^{im\theta})= \left\{
\begin{array}{ll}
0&,~m\neq 0 \\
-i{1+\rho^{2n}\over 1-\rho^{2n}},~m\neq 0,
\end{array}
\right.
\end{equation}
and
\begin{equation}
e^{im\theta}\mapsto {\cal K}_\rho(e^{im\theta})= \left\{
\begin{array}{ll}
0&,~m\neq 0 \\
-2i{\rho^{n}\over 1-\rho^{2n}},~m\neq 0.
\end{array}
\right.
\end{equation}
Furthermore, the radius $\rho$ is explicitly given by
\begin{equation}
\rho=\exp{\left({1\over 2\pi}\log{\left(
\int_{0}^{2\pi}\log{r_1(\phi_1(\theta))} ~d\theta\right)}\right)}.
\end{equation}
\end{thm}
It is straightforward to show that the Theodersen and Garrick
equations enables us to get (\ref{claim2}) when the radius
perturbation $\delta(\theta)$ belongs to $W^{2,\infty}$.

\subsection{Approximation of the Dirichlet-to-Neumann map for perturbed disks.}
Let $D$ be an fixed $\varepsilon-$ perturbation of a disk $B$. We
want to estimate the correction term
$\|\Lambda_D(f)-\Lambda_B(f)\|_{\sH^{-1/2}(S^1)}$ with respect
to the perturbation factor $\varepsilon$. We have

\begin{thm}
\label{theoreme:perturbation:dtn} Let $D$ be a regular domain that
is an $\eps$ perturbation of a disk centered at the origin. Let
$\alpha$ be a real in $(0,1)$. Then there exists a constant $ C>0$
such that
\begin{equation}
\label{perturbation:dtn} \| \DtoN{D}(f)
-\DtoN{B}(f)\|_{\sH^{-1/2}(S^1)} \leq C\eps^{\alpha}\|
f\|_{\sH^{1+\alpha}(S^1)}
\end{equation}
holds for the boundary measurement $f$ belonging to the Sobolev
space $\sH^{1+\alpha}(S^1)$.
\end{thm}
The proof is lengthy and  requires some preliminary results that
we will state and prove before.

\paragraph{Computation of the transplanted Dirichlet-to-Neumann map.}

Our main task is to give the analytic expression of
$\DtoN{B_\rho}^t:\sH^{1/2}(\d\dom) \rightarrow \sH^{-1/2}(\d\dom)$
defined by
\begin{equation*}
\DtoN{B_\rho}^t(F^e) = \sigma_1 \d_r U_{B_\rho}^e,
\end{equation*}
where we set $F^e=f\circ \Psi^e_{|\d\dom}=f\circ\psi^e$. The work
will be divided in two parts: in the first one, we give some
preliminary results based on the expression of
$h=(U_{B_\rho}^e)_{|\d B_\rho}$. This will allow us to get
$\DtoN{B_\rho}^t(F^e)$ and a convenient  approximation
$\DtoN{B_\rho}(f)$. A straightforward calculation shows that for
$\rho<r<1$ we have
\begin{eqnarray*}
U^e(re^{i\theta}) &=& \cfrac{\ln{r}}{\ln{\rho}} \left(c_0(h)-c_0(F^e)
\right) + c_0(F^e)
+ \sum_{n \ne 0} \cfrac{1}{1-\rho^{2|n|}} \left[ r^{|n|} -
\cfrac{\rho^{2|n|}}{r^{|n|}} \right] c_n(F^e) e^{in\theta} \\
&& + \sum_{n \ne 0} \cfrac{\rho^{|n|}}{\rho^{2|n|}-1} \left[
r^{|n|} - \cfrac{1}{r^{|n|}} \right] c_n(h)
e^{in\theta};
\end{eqnarray*}
hence after identification of the Fourier coefficients, it comes that
\begin{equation}\label{fourier_lambda}
\left\{
\begin{array}{ll}
c_0\left(\DtoN{B_\rho}^t(F^e)\right) & = \cfrac{c_0(h)-c_(F^e)}{\ln{\rho}}=0, \\
c_n\left(\DtoN{B_\rho}^t(F^e)\right) & =
\sigma_1\cfrac{|n|}{1-\rho^{2|n|}} \left[ (1+\rho^{2|n|}) c_n(F^e)
-2 \rho^{|n|} c_n(h)\right], ~n\neq 0.
\end{array}
\right.
\end{equation}
We see that the knowledge of $h$ determines uniquely the
Dirichlet-to-Neumann operator $\DtoN{B_\rho}$. It is then useful to get some informations about $h$.
First, we have

\begin{prop}
\label{proposition:equation:h}
We have
\begin{equation}
\label{equation:h}
\begin{split}
\sum_{n \ne 0} |n| c_n(h\circ\xi) e^{in\theta} & + 2 \xi'(\theta) \cfrac{\sigma_1}{\sigma_2} \sum_{n \ne 0}
\cfrac{1+\rho^{2|n|}}{1-\rho^{2|n|}} |n| c_n(h) e^{in\xi(\theta)} \\
&
= 2 \xi'(\theta) \cfrac{\sigma_1}{\sigma_2}
\sum_{n \ne 0} \cfrac{\rho^{|n|}}{1-\rho^{2|n|}} |n| c_n(F^e)
e^{in\xi(\theta)}.
\end{split}
\end{equation}
\end{prop}

\begin{proofof}{Proposition \ref{proposition:equation:h}}
Since $U^i$ is solution of a Dirichlet problem in the disk
$B_\rho$, we obtain
\begin{equation*}
\d_r U^i(\rho e^{i\theta}) = {1\over \rho }\sum_{n \ne 0} |n|
c_n(h\circ \xi) e^{in\theta},
\end{equation*}
and  thanks to the jump
condition satisfied by the normal derivatives $\d_r U^i_{|\d
B_\rho}$ and $\d_r U^e_{|\d B_\rho}$, we deduce equation \eqref{equation:h}.
\end{proofof}
To solve \eqref{equation:h} is a difficult task since the explicit expression of $h$ is hard to manipulate; however when the perturbation factor $\varepsilon$  is very small, one can  give a suitable approximation. Before entering in the details, we have to give some qualitative properties of $h$.

\paragraph{A regularity result.}
We first prove the
following tangential regularity result for solution of the
conductivity problem.

\begin{lemme}
\label{lemme:regularite:tangentielle} Let $u$ be the solution of
the problem P[D,f] \eqref{equation:conductivite} where $\d D$ is
assumed to be  of class $\mathcal{C}^2$.  Then $u_{|\d D}$
belongs to  $\sH^{1+\delta}(\d D,\R)$ for some $\delta \in (0,1)$.
\end{lemme}

\begin{proofof}{Lemma \ref{lemme:regularite:tangentielle}}
From classical methods, the problem $P[D,f]$ has a unique solution
in the variational space $\sH^1(\dom)$ with a trace
$u_{|\d \dom}=f \in \sH^{1/2}(\d \dom)$ and a normal derivative
$\d_n u_{|\d \dom}:=g \in \sH^{-1/2}(\d \dom)$. For all $x \in \d D$, we use the classical
representation formulae for harmonic functions with the help of the
single layer and double layer potential: since  $u$ is harmonic  in $\dom \setminus
D$ and in $D$ we have
\begin{eqnarray*}
\frac{1}{2} u^+(x) &=& \int_{\d\dom} \d_n G(x,y) f(y) ds(y) - \int_{\d D} \d_n G(x,y) u^+(y) ds(y) \\
& & - \int_{\d\dom} G(x,y) g(y) ds(y) +\int_{\d D} G(x,y) \d_n u^+(y) ds(y), \\
\frac{1}{2} u^-(x) &=& \int_{\d D} \d_n G(x,y) u^-(y) ds(y) -\int_{\d D} G(x,y) \d_n
u^-(y)
ds(y)
\end{eqnarray*}
where $G$ is the Newtonian potential and where the normal $n$ to $\partial D$ is oriented to the exterior.
Using the jump conditions $[u]=[\sigma \d_n u]=0$ across the
interface $\d D$, we check that $u=u^+=u^-$ solves the integral
equation
\begin{eqnarray*}
\frac{1}{2} u(x) &+& \cfrac{\sigma_2-\sigma_1}{\sigma_1+\sigma_2}\int_{\d D} \d_n G(x,y) u(y) ds(y) \\
&=& \cfrac{\sigma_1}{\sigma_1+\sigma_2} \left[\int_{\d\dom} \d_n G(x,y) f(y) ds(y) - \int_{\d\dom} G(x,y) g(y)
ds(y)\right].
\end{eqnarray*}

Since $\d D \cap \d \dom = \emptyset$, the right hand side of this
equation is of class $\mathcal{C}^{\infty}$. From the fact that  the boundary $\d D$ is of class $\mathcal{C}^2$, the double layer potential on $\d D$ maps $\sH^{s}(\d D)$ into $\sH^{s+1}(\d D)$ for and hence is compact as operator from $\sH^{s}(\d D)$ into itself.  We conclude thanks to the Fredholm alternative.
\end{proofof}

\paragraph{The perturbation argument.}
Let $T_\xi: \sH^{1+\alpha}(\d B_\rho) \rightarrow \sH^{-1/2}(\d
B_\rho)$ the operator defined by
\begin{equation*}
h\mapsto T_\xi(h) (\theta) = \sum_{n \ne 0} |n|
c_n(h\circ \xi) e^{in\theta} + \xi'(\theta)
\cfrac{\sigma_1}{\sigma_2} \sum_{n \ne 0}
\cfrac{1+\rho^{2|n|}}{1-\rho^{2|n|}} c_n(h) e^{in\xi(\theta)},
\end{equation*}
and $T: \sH^{1+\alpha}(\d B_\rho) \rightarrow \sH^{-1/2}(\d B_\rho)$
the operator defined by
\begin{equation*}
h\mapsto T(h) (\theta) = \sum_{n \ne 0} |n|
c_n(h) e^{in\theta} + \cfrac{\sigma_1}{\sigma_2} \sum_{n \ne 0}
\cfrac{1+\rho^{2|n|}}{1-\rho^{2|n|}} c_n(h) e^{in\theta}.
\end{equation*}
We set $DT_{\xi} =T_\xi -T$. We have

\begin{prop}
\label{proposition:resolubilite:equation:h}
There exists a constant $C>0$ such that
\begin{equation}
\label{majoration:erreur:T} \|T_\xi(u) -T(u) \|_{\sH^{-1/2}(\d
B_\rho)} \leq C \| \xi-I\|_{W^{1}_\infty}\|
u\|_{\sH^{1+\alpha}(\d
B_\rho)}
\end{equation}
holds for all $u$ belonging to $\sH^{1+\alpha}(\d
B_\rho)$ with $0<\alpha<1$.
\end{prop}

\begin{proofof}{Proposition \ref{proposition:resolubilite:equation:h}}
We decompose $ DT_{\xi}(u)= T_1(u)+ T_2(u) + T_3(u)$
where
\begin{eqnarray*}
T_1(u)(\theta) &=& \sum_{n \ne 0} |n| c_n(u\circ
\xi -u) e^{in\theta}, \\
T_2(u)(\theta) &=& -\cfrac{\sigma_1}{\sigma_2}
\left(\xi'(\theta) -1 \right) \sum_{n \ne 0} |n|
\cfrac{1+\rho^{2|n|}}{1-\rho^{2|n|}} c_n(u) e^{in\xi(\theta)},\\
T_3(u)(\theta) &=& -\cfrac{\sigma_1}{\sigma_2} \sum_{n \ne 0}
|n| \cfrac{1+\rho^{2|n|}}{1-\rho^{2|n|}} c_n(u)
\left[e^{in\theta}-e^{in\xi(\theta)}\right].
\end{eqnarray*}

We begin to estimate $\| T_2(u)\|_{\sH^{-1/2}(\d B_\rho)}$, we
have
\begin{equation*}
\| T_2(u)\|_{\sH^{-1/2}(\d B_\rho)} \leq C(\sigma_1,\sigma_2) \|
\xi'-1\|_\infty^2 \|g_\xi\|_{\sH^{-1/2}(\d B_\rho)},
\end{equation*}
where $g_\xi = g_1\circ \xi + g_2\circ \xi$ with
\begin{equation*}
g_1(\theta) = \sum_{n\ne 0} |n| c_n(u) e^{in\theta}
\end{equation*}
and
\begin{equation*}
g_2(\theta) = 2\sum_{n\ne 0} |n| \cfrac{\rho^{2|n|}}{1-\rho^{2|n|}} c_n(u)
e^{in\theta}.
\end{equation*}
While the estimation of $\|g_2\circ\xi\|_{\sH^{-1/2}(\d B_\rho)}$
is straightforward
\begin{equation*}
\|g_2\circ\xi\|_{\sH^{-1/2}(\d B_\rho)} \leq C(\delta_0) \|u\|_{\sH^{1/2}(\d
B_\rho)},
\end{equation*}
the estimation of $\|g_2\circ\xi\|_{\sH^{-1/2}(\d B_\rho)}$ is a
little bit harder. We first observe that, since $u$ belongs to
$\sH^s(\d\dom), ~s>1$, we can define $\mathcal{H} u' =
(\mathcal{H}s)'$ where $\mathcal{H}$ is the Hilbert transform on
the circle. It then comes that
\begin{equation*}
g_1\circ \xi = \mathcal{H}u'\circ \xi = (\mathcal{H}u)'\circ
\xi,
\end{equation*}
and then that
\begin{equation*}\begin{split}
\| g_1\circ \xi \|_{\sH^{-1/2}(\d B_\rho)} \leq C
\|(\mathcal{H}u)'\circ \xi\|_{\sH^{-1/2}(\d B_\rho)} &\leq C
\|(\mathcal{H}u)'\|_{\sH^{-1/2}(\d B_\rho)} \\
&\leq C \|\mathcal{H}u\|_{\sH^{1/2}(\d B_\rho)} \leq C
\|u\|_{\sH^{1/2}(\d B_\rho)}.
\end{split}
\end{equation*}
In the last inequality, we used the fact that the Hilbert
transform corresponds to an unimodular multiplier and then is an
isometry on $\sH^{1/2}$. Gathering the estimates on $g_1\circ\xi$
and $g_1\circ\xi$, we get
\begin{equation*}
\| T_2(u)\|_{\sH^{-1/2}(\d B_\rho)} \leq
C(\sigma_1,\sigma_2,\delta_0) \| \xi'-1\|_\infty^2
\|u\|_{\sH^{1/2}(\d B_\rho)},
\end{equation*}
It remains to estimate $\| T_3(u)\|_{\sH^{-1/2}(\d B_\rho)}$, we
have
\begin{equation*}
\begin{split}
\| T_3(u)\|_{\sH^{-1/2}(\d B_\rho)} \leq C(\sigma_1,\sigma_2)&
\| \sumZe |n| c_n(u) \left( e^{in\xi(\theta)} -e^{in \theta} \right) \|_{\sH^{-1/2}(\d B_\rho)} \\
& +
\| \sumZe |n| \cfrac{\rho^{2|n|}} { 1-\rho^{2|n|} }c_n(u) \left(
e^{in\xi(\theta)} -e^{in \theta} \right) \|_{\sH^{-1/2}(\d
B_\rho)}.
\end{split}\end{equation*}
We focus on the first part of the sum. We have
\begin{eqnarray*}
\sumZe |n| c_n(u) \left( e^{in\xi(\theta)} -e^{in \theta} \right) &=& (\mathcal{H}u)' \circ \xi -\mathcal{H} u' \\
&=& (\mathcal{H}u' \circ \xi)\xi' + (1-\xi') (\mathcal{H} u')\circ\xi -(\mathcal{H}u)' \circ \xi - (\mathcal{H} u)' \\
&=& \left( (\mathcal{H}u) \circ \xi - \mathcal{H} u \right)' +
(1-\xi') (\mathcal{H} u')\circ\xi.
\end{eqnarray*}
Hence,
\begin{equation*}
\begin{split}
\| \sumZe |n| c_n(u) &( e^{in\xi(\theta)} - e^{in \theta} )\|_{\sH^{-1/2}(\d B_\rho)} \\
&\leq \| (\mathcal{H}u) \circ \xi - \mathcal{H} u \|_{\sH^{1/2}(\d B_\rho)}+ \|\xi'-1\|_{\infty} \|u'\|_{\sH^{-1/2}(\d B_\rho)} \\
&\leq \| \sumZe c_n(\mathcal{H}u) ( e^{in\xi(\theta)} -e^{in
\theta} ) \|_{\sH^{1/2}(\d B_\rho)}+ \|\xi'-1\|_{\infty} \|u
\|_{\sH^{1/2}(\d B_\rho)}
\end{split}
\end{equation*}
Suppose an instant that $u$ is a trigonometric polynomial of degree
$d$. Thanks to the composition Theorem
\ref{theoreme:precomposition}, for $\delta \in(0,1)$, we have
\begin{equation*}
\begin{split}
\| \sumZe c_n(\mathcal{H}u) ( e^{in\xi(\theta)} -e^{in \theta} ) \|_{\sH^{-1/2}(\d B_\rho)} &\leq \sum_{|n|\leq d} |c_n(u)|
\| e^{in\xi(\theta)} -e^{in \theta} \|_{\sH^{1/2}(\d B_\rho)} \\
&\leq  C(\delta)~ \|u\|_{\sH^{1+\delta}(\d B_\rho)} ~
\omega_{\delta}(\|\phi-I\|_{\sW^{1,\infty}(S^1)}),
\end{split}
\end{equation*}
where $\omega_{2\delta}$ is the modulus of continuity defined by
\eqref{definition:omega}. Then, there exists a constant
$C(\sigma_1,\sigma_2,\delta_0)$ such that
\begin{equation*}
\| T_3(u)\|_{\sH^{-1/2}(\d B_\rho)} \leq C \omega_\delta
(\|\xi-Id\|_{\sW^{1,\infty}(\d B_\rho)}) \|u\|_{\sH^{1+\delta}(\d
B_\rho)}.
\end{equation*}
The estimation of $T_1(u)$ obeys to the same computation : we have
by following the same lines
\begin{equation*}
\| T_1(u)\|_{\sH^{-1/2}(\d B_\rho)} \leq C \omega_\delta
(\|\xi-Id\|_{\sW^{1,\infty}(\d B_\rho)}) \|u\|_{\sH^{1+\delta}(\d
B_\rho)}
\end{equation*}
end this ends our proof.
\end{proofof}

\begin{proofof}{Theorem \ref{theoreme:perturbation:dtn}}
Let us return to the equations satisfied by $h$. Set
\begin{equation*}
b(F^e,\rho,\xi,\sigma_1,\sigma_2)(\theta) = 2\xi'(\theta)
\cfrac{\sigma_1}{\sigma_2} \sumZe \cfrac{\rho^{|n|}}{1-\rho^{2|n|}
} |n| c_n(F^e) e^{in\xi(\theta)},
\end{equation*}
 and suppose that the perturbation parameter $\varepsilon$ is
sufficiently small such that
\begin{equation*}
\| (DT_\xi)T^{-1}\|<1.
\end{equation*}
It follows that
\begin{equation*}
\begin{array}{lll}
h (\theta)&=& T_\xi^{-1}b(F^e,\rho,\xi,\sigma_1,\sigma_2)\\
&=&T^{-1}\left( I+DT_{\xi}T^{-1}
\right)^{-1}b(F^e,\rho,\xi,\sigma_1,\sigma_2)\\
&=&T^{-1} b(F^e,\rho,\xi,\sigma_1,\sigma_2)-T^{-1}\left(
I+DT_{\xi}T^{-1}
\right)^{-1}DT_{\xi}T^{-1}b(F^e,\rho,\xi,\sigma_1,\sigma_2).
\end{array}
\end{equation*}
An easy calculation of $T^{-1}$ shows that
\begin{equation*}
h(\theta) = 2 \cfrac{\sigma_1}{\sigma_2} \sumZe
{\cfrac{\rho^{|n|}}{1-\rho^{2|n|} } \over \displaystyle\rho^{|
n|}+2{\sigma_1\over \sigma_2}{1+\rho^{2|n|}\over
1-\rho^{2|n|}}} c_n(f) e^{in \theta} +\delta u(\theta),
\end{equation*}
where $\delta u\in\mathcal{C}^{\infty}$ is such that for all $s>0$
$\|\delta u\|_{\sH^m(S^1)} \leq C \| \xi-I\|_{\sW^1_\infty(S^1)}\| f
\|_{\sH^{1/2}}$ for all integers  $m \in \mathbb{N}$.
Plugging this expression of $h$ in  formula (\ref{fourier_lambda})
and setting $\mu={\sigma_2-\sigma_1\over \sigma_1+\sigma_1}$ it
follows  that
\begin{equation}\label{approx_DtN}
\left\{
\begin{array}{ll}
c_0\left(\DtoN{B_\rho}^t(F^e)\right) & = \cfrac{c_0(h)-c_(F^e)}{\ln{\rho}}=0, \\
c_n\left(\DtoN{B_\rho}^t(F^e)\right) & =
|n|\sigma_1\left(\displaystyle{1+\mu \rho^{2|n|}\over
1-\mu\rho^{2|n|}}c_n(f)+~c_n(\delta h)\right)n\neq 0.
\end{array}
\right.
\end{equation}
where $\delta h$ is a function such that $\| \delta
h\|_{H^1/2(S^1)} \le C \varepsilon$.

We are ready now to finish the proof. First of all, we know from
the Cauchy-Riemann identities that
\begin{equation*}
\DtoN{B_\rho}^t(F^e)=\left(g\circ\psi^e\right)(\psi^e)'
\end{equation*}
and using the same arguments that we developed above, one can
easily show that there exists a positive constant $C>0$ depending
on $\rho_0$ and $\delta_0$ such that
\begin{equation}
\| \DtoN{B_\rho}^t(F^e)-\Lambda_D(f)\|_{\sH^{-1/2}(S^)}
\le C \varepsilon^\alpha \| f \|_{\sH^{1+\alpha}}.
\end{equation}
Hence we get
\begin{equation*}
\begin{array}{lll}
\| \DtoN{B}(f)-\DtoN{D}(f)\|_{\sH^{-1/2}(S^1)} &=&
\| \DtoN{B}(f)-\DtoN{B_\rho}^t(F^e)+\Lambda_{B_\rho}^t(f)-\Lambda_D(f)\|_{\sH^{-1/2}(S^1)},\\
&\le&\| \DtoN{B}(f)-\DtoN{B_\rho}^t(F^e)+ C\varepsilon^\alpha \| f \|_{\sH^{1+\alpha}(S^1)},\\
&\le &\|\DtoN{B}(f)-\DtoN{B_\rho}(f)\|_{\sH^{-1/2}(S^1)}+C \|\left({\cal H}(\delta
h)\right)'\|_{\sH^{-1/2}(S^1)},\\
&&~~~+C\varepsilon^\alpha \| f \|_{\sH^{1+\alpha}(S^1)},\\
&\le
&\|\DtoN{B}(f)-\DtoN{B_\rho}(f)\|_{\sH^{-1/2}(S^1)}+C\varepsilon^\alpha
\| f \|_{\sH^{1+\alpha}(S^1)}.
\end{array}
\end{equation*}
Recall that the Fourier coefficients of $\DtoN{B_\rho}(f)$ are
given by 
\begin{equation*}
c_n(\DtoN{B_\rho}(f))= |n|\sigma_1 \displaystyle{1+\mu
\rho^{2|n|}\over 1-\mu\rho^{2|n|}}c_n(f).
\end{equation*}
 Hence after denoting
$\rho_1$ the radius of the disk $B$, we get
\begin{equation}
\Lambda_{D}(f)-\Lambda_{B_\rho}(f)=\sigma_1\sum_{n\neq
0}|n|\left({1+\mu \rho^{2|n|}\over 1-\mu \rho^{2| n|}} -{1+\mu
\rho_1^{2|n|}\over 1-\mu \rho_1^{2| n|}} \right)c_n(f)
e^{in\theta}
\end{equation}
and this implies that we can find  a constant $C>0$ depending on
$\delta_0$,$\rho_0$ and on the conductivities $\sigma_i,~i=1,2$
such that
\begin{equation}
\|\Lambda_{D}(f)-\Lambda_{B_\rho}(f)\|_{\sH^{-1/2}(S^1)} \le C|
\rho-\rho_1|\|  f \|_{\sH^{1/2}(S^1)};
\end{equation}
 we conclude thanks to the fact that $| \rho_1-\rho_2| \le
C\varepsilon$.  This ends the proof of our theorem.
\end{proofof}

%%%%%%%%%%%%%%%%%%%%%%%%%%%%%%%%%%%%%
\section{Proof of the main theorems}
%%%%%%%%%%%%%%%%%%%%%%%%%%%%%%%%%%%%%

We subdivide  the section in two parts : in the first one, we
focus on the case where the inclusions are disks. In the second
part, the inclusions belong to $C[ \varepsilon]$ .

\subsection{Identifiability for disks : proof of Theorem \ref{unicite:disques}}

For the case where the domains are disks, we can always suppose
that $D_1$ is centered at the origin.
\begin{proofof}{Theorem \ref{unicite:disques}}
We use conformal mappings that maps a non concentric disk $D_2$
into a disk centered at the origin. We know that this can be done
by means of the Moebius transform
\begin{equation*}
w(z) =\cfrac{z-b}{1-\overline{b}z} \text{ with } |b|<1.
\end{equation*}
In \cite{Ahlfors,Henrici1,Henrici3}, the interested
reader will find all the details about the properties on such
transforms . The radius of the transformed disk is denoted by
$R_2$. If $z=re^{i\theta}$, then $w(z)$ writes
$\rho(r,\theta)e^{i\phi(r,\theta)}$ where
\begin{equation*}
\rho^2(r,\theta) = \cfrac{r^2+|b|^2-r(\overline{b} e^{i\theta} +b e^{-i\theta}) }{1+|b|^2-r(\overline{b} e^{i\theta} +b e^{-i\theta})},
\end{equation*}
and where
\begin{equation*}
\phi(r,\theta) = \arctan{\cfrac{r\sin{\theta} -\Im{b} r^2+ r
[\sin{\theta} (\Im{b}^2-\Re{b}^2)+2\Re{b}\Im{b}\cos{\theta}]}
{2\cos{\theta}-\Re{b}r^2-\Re{b}+r[(\Re{b}^2-\Im{b}^2) \cos{\theta}
-2 \Re{b}\Im{b}\sin{\theta} ]}}.
\end{equation*}
Here $\Re{b}$ and $\Im{b}$ denote the real and imaginary part of the complex $b$.
A straightforward computation shows that
\begin{equation*}
\d_r \rho(1,\theta) =
\cfrac{1-|b|^2}{1+|b|^2-(\overline{b}e^{i\theta} + be^{-i\theta}
)} \text{ and } \d_r \phi(1,\theta) = 0.
\end{equation*}
From the chain rule of differentiation, we get:
\begin{equation*}
\Lambda_{D_2}(f) = \d_r \rho(1,\theta) \Lambda_C^{R_2}(f\circ \phi^{-1}),
\end{equation*}
where $\Lambda_C^{R_2}$ is the Dirichlet-to-Neumann map for the
transformed and concentric problem. It follows that
\begin{equation*}
\Lambda_{D_1}(f) =\Lambda_{D_2}(f)\Rightarrow \Lambda_{D_1}(f) =
\cfrac{1-|b|^2}{1+|b|^2-(\overline{b}e^{i\theta} + b ^{-i\theta}
)} \Lambda_C^{R_2}(f\circ \phi^{-1})
\end{equation*}
and the key for solving the problem is to use the Dirichlet-to-Neumann map for concentric disk. We have:
\begin{equation*}
\cfrac{1+\mu R_1^2}{1-\mu R_1^2} \cos{\theta} = \cfrac{1-|b|^2}{1+|b|^2-(\overline{b}e^{i\theta} + b ^{-i\theta} )} \sum_{k \ne 0} |k| \cfrac{1+\mu R_2^{2|k|}}{1-\mu R_2^{2|k|}} c_k(f\circ \phi^{-1}) e^{ik\phi(\theta)},
\end{equation*}
or equivalently:
\begin{equation*}
\begin{split}
\cfrac{1+\mu R_1^2}{1-\mu R_1^2} &\left[ (1+|b|^2) \cos(\theta) - \Re(b) - \frac{1}{2} (\overline{b} e^{2i\theta}+be^{-2i\theta}) \right] \\
&= (1-|b|^2) \sum_{k \ne 0} |k| \cfrac{1+\mu R_2^{2|k|}}{1-\mu R_2^{2|k|}} c_k(f\circ \phi^{-1}) e^{ik\phi(\theta)}
\end{split}
\end{equation*}
Replacing $\theta $ by $\phi^{-1}(\theta)$, we then get that
$\Lambda_{D_1}(f)=\Lambda_{D_2}(f)$ implies that the $2\pi$
periodic function
\begin{equation*}
\begin{split}
F(\theta) = (1-|b|^2) \sum_{k \ne 0} |k| &\cfrac{1+\mu
R_2^{2|k|}}{1-\mu R_2^{2|k|}} c_k(f\circ \phi^{-1})
e^{ik\theta}-\cfrac{1+\mu R_1^2}{1-\mu R_1^2} \\&\left[ (1+|b|^2)
\cos(\phi^{-1}(\theta)) - \Re(b) - \frac{1}{2} (\overline{b}
e^{2i\phi^{-1}(\theta)}+be^{-2i\phi^{-1}(\theta)}) \right]
\end{split}
\end{equation*}
satisfies $c_k(F)=0$ for all $k\in\Z$.
\par
\noindent We tackle the computation of these  Fourier
coefficients. First of all, we need to compute
$c_k(f\circ\phi^{-1})$ and $c_k(e^{2i\phi^{-1}})$. A
straightforward computation shows that:
\begin{equation*}
c_k(f\circ\phi^{-1}) =
\left\{
\begin{array}{l}
\Re(b)\text{ if }k=0,\\
\cfrac{1}{2}(1-|b|^2)(-b)^{k-1}\text{ else} ;
\end{array}
\right.
\end{equation*}
and that
\begin{equation*}
c_k(e^{i\phi^{-1}}) = \left\{
\begin{array}{l}
b\text{ if }k=0,\\
(1-|b|^2)(-\overline{b})^{k-1}\text{ if } k>0,\\
0 \text{ else}.
\end{array}\right.
\end{equation*}
Since
\begin{equation*}
e^{i\phi^{-1}(\theta)} = b+ (1-|b|^2) \cfrac{e^{i\theta}}{1+\overline{b} e^{i\theta}},
\end{equation*}
we deduce that
\begin{equation*}
e^{2i\phi^{-1}(\theta)} = b^2 + 2b (1-|b|^2) \cfrac{e^{i\theta}}{1+\overline{b} e^{i\theta}} + (1-|b|^2)^2 \cfrac{e^{2i\theta}}{(1+\overline{b} e^{i\theta})^2},
\end{equation*}
and then
\begin{equation*}
c_k(e^{2i\phi^{-1}(\theta)}) = \left\{
\begin{array}{l}
b^2 \text{ if }k=0,\\
(-\overline{b})^{k-2} (1-|b|^2) \left[k-1 -(k+1)|b|^2\right]\text{ if } k>0,\\
0 \text{ else}.
\end{array}\right.
\end{equation*}
Hence, the Fourier coefficients $(c_k(F))_k$ are given by the
following formulae
\begin{equation*}
c_0(F) =0 \text{ and } c_1(F) = \mu (1-|b|^2)^2 \cfrac{R_2^2-R_1^2}{(1-\mu R_2^2)(1-\mu R_1^2)}
\end{equation*}
and for $k\geq 2$:
\begin{equation*}
\begin{split}
&c_k(F) = \cfrac{(-1)^{k-1}(1-|b|^2)}{2} \Big[ k \cfrac{1+\mu R_2^{2|k|}}{1-\mu R_2^{2|k|}} (1-|b|^2) b^{k-1} - \cfrac{1+\mu R_1^{2}}{1-\mu R_1^{2}} (1+|b|^2) b^{k-1}\\
&~~~~~~~~~~~~~~~~~~~~~~~~~~~~~~~~~~~~~~- \cfrac{1+\mu R_1^{2}}{1-\mu R_1^{2}}\left[k-1-|b|^2 (k+1)\right] \overline{b}^{k-1}\Big],\\
&= \cfrac{(-1)^{k-1}(1-|b|^2)}{2} \left[ b^{k-1} 2\mu k
\cfrac{(R_2^{2|k|} -R_1^2)(1-|b|^2)}{(1-\mu R_2^{2|k|})(1-\mu
R_1^2)} + (\overline{b}^{k-1}-b^{k-1}) (1+|b|^{2}) \cfrac{1+\mu
R_1^2}{1-\mu R_1^2}\right].
\end{split}
\end{equation*}
Let us show that $c_k(F)=0,~\forall k \ge 0$ implies $D_1=D_2$.
First of all, the condition $c_1(F)=0$ implies that $R_1=R_2$. It
remains to show that the Moebius transform is, in fact, the
identity. Equivalently, we need to prove that $b=0$.

Let us assume that we have $b\ne 0$; the condition
$c_k(F)=0,~\forall k \ge 2 $ would imply that $\overline{b}/b$ is
real and then that $\overline{b} =\pm b$. From the identity
\begin{equation*}
2\mu \cfrac{(R_1^{2|k|}) (1-|b|^2)}{(1-\mu R_2^{2|k|})(1-\mu R_1^2)} + \left[\left(\cfrac{\overline{b}}{b}\right)^{k-1}-1\right] (1+|b|^{2}) \cfrac{1+\mu R_1^2}{1-\mu R_1^2}=0,
\end{equation*}
we would get for $k=2p+1, ~p\ge 1$
\begin{equation*}
\cfrac{2\mu(R_2^{4p+2}-R_1^2)(1-|b|^2)}{(1-\mu R_2^{4p+2})(1-\mu
R_1^2)} =0
\end{equation*}
and then  $| b| =1$;  this  is impossible since the Moebius
transform requires $|b|<1$. Hence $b=0$ is the only possibility.
Gathering the two identities $R_1=R_2$ and $b=0$, it then comes
that $D_1=D_2$.
\end{proofof}

\subsection{Proof of theorem \ref{unicite:fabes} : the approximate identifiability for approximate disks}

In this section, we need an intermediary lemma showing  the
relations between the norm of the superposition operators
generated by a diffeomorphism $\xi$ and its inverse $\xi^{-1}$.
Its proof can found in (\cite{Bourdaud}); we have

\begin{lemme}
\label{lemme:bourdeau} Let $T_\xi: \sH^{1/2}(S^1) \rightarrow
\sH^{1/2}(S^1)$ be the composition operator defined by
$T_\xi(f)=f\circ\xi$. Then, if $T_\xi$ is bounded on
$\sH^{1/2}(S^1)$, then so is $T_{\xi^{-1}}$. Furthermore, we have
\begin{equation*}
\|T_\xi\|_{\mathcal{L}(\sH^{1/2}(S^1),\sH^{1/2}(S^1)} =
\|T_{\xi^{-1}}\|_{\mathcal{L}(\sH^{1/2}(S^1),\sH^{1/2}(S^1)}.
\end{equation*}
\end{lemme}

\begin{proofof}{Theorem \ref{unicite:fabes}} If $\Psi_1$ is the conformal mapping
which maps $\dom\setminus \overline{D}_1$ onto an annulus
$A(1,R_1^{\eps}) = \dom \setminus B(0,R_1^\eps)$, it then comes from
the results of the preceding section that
\begin{equation*}
\| \Lambda_{B^\eps_1}(f\circ \varphi_{1,\eps}) - \Lambda_{D_1}(f)
\|_{\sH^{-1/2}} \leq C \eps^\alpha, ~0<\alpha<1
\end{equation*}
where $\varphi_{1,\eps} = \Psi_1^{-1}~_{|\d\dom}$ and where
$B^\eps_1$ denotes the disk $B(0,R_1^\eps)$.

Concerning $D_2$ which is the $\eps$-perturbation of a non
concentric disk $B_2$, we begin to transform it via $\Psi_2$ the
Moebius transform that maps $B_2$ onto a concentric disk of radius
$R_2$. It is obvious that $\tilde{D_2}=\Psi_2(D_2)$ is a slight
perturbation of the concentric disk $\tilde{B_2}= \Psi_2(B_2)$.

The chain rule derivative gives
\begin{equation*}
\Lambda_{D_2}(f) = \cfrac{1-|b_2|^2}{1+|b_2|^2-\overline{b}_2
e^{i\theta} -b_2 e^{-i\theta}}
\Lambda_{\tilde{D_2}}(f\circ\varphi_2),
\end{equation*}
where $\varphi_2 = \Psi_2^{-1}~_{|\d\dom}$. Thanks to our
preceding results, we know that if $\Omega\backslash\overline{
B(0,R_2^\eps)}$ is the conformal transform of
$\Omega\backslash\overline{\tilde{D_2}}$ then
\begin{equation*}
\| \Lambda_{B(0,R_2^\eps)}(f\circ \varphi_{1,\eps}) -
\Lambda_{\tilde{D_2}}(f) \|_{\sH^{-1/2}} \leq C \eps^\alpha,
~0<\alpha<1
\end{equation*}
and the assumption $\DtoN{D_1}(f)=\DtoN{D_2}(f)$ implies that
\begin{equation*}
\begin{split}
\| \DtoN{B(0,R_1^\eps)}(f) &- \cfrac{1-|b_2|^2}{1+|b_2|^2-\overline{b_2} e^{i\theta} -b_2 e^{-i\theta} } \DtoN{B(0,R_2^\eps)}(f\circ \phi_2) \|_{\sH^{-1/2}(S^1)} \\
& \leq \| \DtoN{D_1}(f) - \DtoN{B(0,R_1^\eps)}(f)\|_{\sH{-1/2}(S^1)} \\
&+ \sup_{\theta \in [0,2\pi]} \left| \cfrac{1-|b_2|^2}{1+|b_2|^2-\overline{b_2} e^{i\theta} -b_2 e^{-i\theta} }\right| \| \DtoN{\tilde{D}_2}(f\circ\phi_2 ) - \DtoN{B(0,R_2^\eps)}(f\circ\phi_2)\|_{\sH{-1/2}(S^1)} \\
& \leq C \eps^\alpha
\end{split}
\end{equation*}
where $C=C(\delta_0)$ is a constant depending only on $\delta_0$
and where $0<\alpha<1$. Following the same lines for the proof of
Theorem \ref{unicite:disques}, we show that
$\Lambda_{D_1}(f)=\Lambda_{D_2}(f)$ implies that
\begin{equation*}
\|G\|_{\sH^{-1/2}(S^1)} \leq C\eps^\alpha
\end{equation*}
where $G$ is the $2\pi-$ periodic function defined by
\begin{equation*}
\begin{split}
G(\theta) = &(1-|b_2|^2) \sumZe |k| \cfrac{1+ \mu (R_2^\eps)^{2|k|}}{1- \mu (R_2^\eps)^{2|k|}} c_k(f\circ \phi^{-1}) e^{ik\phi(\theta)}\\
& -\cfrac{1+\mu (R_1^\varepsilon)^2}{1-\mu (R_1^\varepsilon)^2}
\left( (1+|b|^2)\cos{\theta} -\frac{1}{2} \overline{b_2}
e^{2i\theta} - \frac{1}{2} b_2 e^{-2i\theta}-\Re(b_2) \right).
\end{split}
\end{equation*}
We set $F=G\circ\phi^{-1}$; from Lemma \ref{lemme:bourdeau}, we
deduce that
\begin{equation*}
\Lambda_{D_1}(f)=\Lambda_{D_2}(f) \implies \|F\|_{\sH^{-1/2}(S^1)}
\leq C\eps^{\alpha},
\end{equation*}
or equivalently
\begin{equation}
\label{condition:petitesse:somme} \sumZe \cfrac{|c_k(F)|^2}{|k|}
\leq C\eps^{2\alpha}.
\end{equation}
We claim that this last inequality implies the two following
inequalities
\begin{equation*}
|b_2|^2\leq C\eps^2 \text{ and }
|R^\varepsilon_1-R^\varepsilon_2|\leq C\eps^{2\alpha}.
\end{equation*}
Indeed, we write
\begin{equation*}
c_k(F)=\alpha_k+\beta_k+\gamma_k,
\end{equation*}
where for $k\geq 2$
\begin{eqnarray*}
\alpha_k&=& (-1)^{k-1} (1-|b_2|^2)^2 \mu b^{k-1} k \cfrac{(R_2^\varepsilon)^{2k}-(R_1\varepsilon)^{2k}}{(1-\mu(R_1^\eps)^{2}) (1-\mu(R_2^\eps)^{2k}) } ,\\
\beta_k &=& \overline{b}^{k-1} (1+|b_2|^2)^2
\cfrac{1+\mu(R_1^\eps)^{2}}{1-\mu(R_1^\eps)^{2}},
\end{eqnarray*}
and where $\gamma_k = - \overline{\beta}_k$. When $k=1$, we have
\begin{equation*}
c_1(F) = \mu (1-|b_2|^2)^2
\cfrac{(R_2\varepsilon)^2-(R_1^\varepsilon)^2}{(1-\mu R_1^2)
(1-\mu R_2^2)}.
\end{equation*}
From \eqref{condition:petitesse:somme}, we get
\begin{equation*}
|c_1(F)| \leq C \eps^\alpha \text{ and } \sum_{|k|\geq2}
\cfrac{|c_k(f)|^2}{|k|} \leq C\eps^{2\alpha}.
\end{equation*}
As we got for the case of disks, the condition $|c_1(F)| \leq
C\eps^\alpha$ implies
\begin{equation*}
|(R_1^\varepsilon)^2-(R_2^\varepsilon)^2|(1-|b_2|^2)^2\leq C
\eps^{2\alpha}
\end{equation*}
where the constant $C$ depends on $\delta_0$. Since $(1-|b_2|^2)^2
> \delta_1$, we then get $|R_2^\varepsilon-R_1^\varepsilon|\le C \varepsilon^\alpha$
This means that the radii of the two disks are very close.

Let us prove that
\begin{equation*}
\sum_{|k|\geq2} \cfrac{|c_k(f)|^2}{|k|} \leq C\eps^2 \implies |b_2| \leq C
\eps^{2\alpha}.
\end{equation*}
Since
\begin{equation*}
c_k(F) = (-1)^k b_2^{k-1} k \left[ (1-|b_2|^2)^2
\cfrac{(R_2^\varepsilon)^{2k}
-(R_1^\varepsilon)^2}{(1-\mu(R_1^\eps)^{2})
(1-\mu(R_2^\eps)^{2k})} + \left(\left(\cfrac{\overline{b_2}}{b_2}
\right)^{k-1}-1\right) \cfrac{1+\mu (R_1^\varepsilon)^2}{k(1-\mu
R_1^2)} \right]
\end{equation*}
we also get
\begin{equation}\label{delta}
\sum_{k\geq 2} k |b_2|^{2(k-1)} \delta_k \leq \eps^{2\alpha},
\end{equation}
where we set
\begin{equation*}
\delta_k = \left| (1-|b_2|^2)^2 \cfrac{(R_2^\varepsilon)^{2k}
-(R_1^\varepsilon)^2}{(1-\mu(R_1^\eps)^{2})
(1-\mu(R_2^\eps)^{2k})} + \left(\left(\cfrac{\overline{b_2}}{b_2}
\right)^{k-1}-1\right) \cfrac{1+\mu (R_1^\varepsilon)^2}{k(1-\mu
(R_1^\varepsilon)^2)} \right|.
\end{equation*}
Let us check that $\delta_2>0$. We argue by contradiction and
assume the converse: if $\delta_2=0$ then we would have
\begin{equation*}
(1-|b_2|^2)^2 \cfrac{(R_2^\varepsilon)^{4}
-(R_1^\varepsilon)^2}{(1-\mu(R_1^\eps)^{2}) (1-\mu(R_2^\eps)^{4})}
+ \left(\cfrac{\overline{b_2}}{b_2} -1\right) \cfrac{1+\mu
(R_1^\varepsilon)^2}{k(1-\mu (R_1^\varepsilon)^2)} =0.
\end{equation*}
Taking the imaginary part of $b_2$, we would get $\overline{b_2}
=\pm (b_2)^{-1}$ and then
$(R_2^\varepsilon)^4=(R_1^\varepsilon)^2$ which is impossible
since we have $|R_1^\varepsilon-R_2^\varepsilon|\leq
C\eps^{2\alpha}$ with $R_1^\varepsilon,R_2^\varepsilon>\rho_0$.
Hence, from (\ref{delta}), we deduce
\begin{equation}
| b_2|^2 \le C \varepsilon^{2\alpha}
\end{equation}
where $C>0$ is a positive constant that depends on $\delta_0$ and
$\rho_0$.
\par
\noindent Let us sum up our conditions :
\begin{itemize}
\item we have $| (R_1^\varepsilon)^2 -(R^2_\varepsilon|) \le C
\varepsilon^\alpha$; this means that
\begin{equation*}\label{difference_12}
\big| | B_1^\varepsilon | - | B_2^\varepsilon t| \big|\le C
\varepsilon^{2\alpha}
\end{equation*}
with a constant $C>0$ depending on $\delta_0$ and $\rho_0$. \item
We have $| b_2| \le C \varepsilon^\alpha$, this means that
the center $b_2$ of $B_2$ is near the origin $0$; a
straightforward calculus gives
\begin{equation*}\label{difference_eps2}
| B_2\Delta B_2^\varepsilon | \le C \varepsilon^\alpha.
\end{equation*}
\end{itemize}
\par
\noindent
 Since $| D_2\Delta B_2| \le C \varepsilon$, we then get
\begin{equation*}
| D_2\Delta B_2^\varepsilon| \le \varepsilon^\alpha
\end{equation*}

It then follows that
\begin{equation*}
| D_2\Delta B_1^\varepsilon| \le \varepsilon^\alpha
\end{equation*}
and then that
\begin{equation*}
| D_2 \Delta D_1|\le \varepsilon^\alpha.
\end{equation*}
This ends the proof of the result.
\end{proofof}

%%%%%%%%%%%%%%%%%%%%%%%%%%%%%%%%%%
\section{Proof of the precomposition theorem.}
%%%%%%%%%%%%%%%%%%%%%%%%%%%%%%%%%%

The main tool for proving Theorem \ref{theoreme:precomposition} is
the following lemma. It provides the behaviour of the precomposition on the Fourier basis functions.

\begin{lemme}
\label{lemme:composition:atomique} Assume $\phi$ is a
$\sW^{1,\infty}$ diffeomorphism on $S^1$. Then, for all $\delta \in
(0,1/2)$
\begin{equation}\label{inegalite:composition:atomique}
\|e^{in\phi(\theta)} - e^{in\theta} \|_{\sH^{1/2}(S^1)} \leq
C(\delta)~ n^{1+2\delta}~ \omega_{2\delta}(\|\phi-I\|_{\sW^{1,\infty}(S^1)})
\end{equation}
where $\omega_{2\delta}$ is the modulus of continuity defined by
\eqref{definition:omega}.
\end{lemme}

\begin{proofof}{Lemma \ref{lemme:composition:atomique}}
The intrinsic definition of the $\sH^{1/2}(S^1)$ gives
\begin{equation*}
\|e^{in\phi(\theta)} - e^{in\theta} \|_{\sH^{1/2}(S^1)}
=\cfrac{1}{4\pi} \iint_{S^1\times S^1} \cfrac{|e^{in \phi(\theta)}-e^{in\theta} -e^{in
\phi(\alpha)}+e^{in\alpha}|^2}{|e^{i\theta}-e^{i\alpha}|^2}d\theta
d\alpha.
\end{equation*}
Noting that
\begin{equation*}
|e^{in \phi(\theta)}-e^{in\theta} -e^{in \phi(\alpha)}+e^{in\alpha}|^2
\leq 2 \big[ |e^{in(\phi-Id)(\theta)} - e^{in(\phi-Id)(\alpha)}|^2
+|e^{in(\phi-Id)(\alpha)} - 1|^2 |e^{in\theta}-e^{in\alpha}|^2
\big]
\end{equation*}
we write $\|e^{in\phi(\theta} - e^{in\theta} \|_{\sH^{1/2}(S^1)}
\leq I_1+I_2$ where
\begin{equation*}
I_1 =\cfrac{1}{2\pi} \iint_{S^1\times S^1} |e^{in(\phi-Id)(\alpha)}-1|^2\cfrac{|e^{in\theta}-e^{in\alpha}|^2}{|e^{i\theta}-e^{i\alpha}|^2}d\theta
d\alpha,
\end{equation*}
and
\begin{equation*}
I_2 =\cfrac{1}{2\pi} \iint_{S^1\times S^1} \cfrac{|e^{in(\phi-Id)(\theta)}-e^{in(\phi-Id)(\alpha)}|^2}{|e^{i\theta}-e^{i\alpha}|^2}d\theta
d\alpha.
\end{equation*}
Concerning $I_1$, we follow an idea of Bourgain, Brezis and
Mironescu (\cite{BouBreMir}) and first integrate with respect to $\theta$.
Thanks to Parseval's formula, we get for all $\delta \in (0,1)$
\begin{equation*}
I_1 = n \int_{S^1} |e^{in(\phi-Id)(\alpha)}-1|^2 d\alpha \leq C
~n^{1+2\delta}~ \|\phi-Id\|_{\infty}^{2\delta}.
\end{equation*}
To estimate $I_2$, we introduce a non-negative parameter $l$ and
split $I_2$ into $I_2^d+I_2^r$ where
\begin{equation*}
I_2^d = \cfrac{1}{2\pi} \iint_{|\theta-\alpha|<l}\left(
\cfrac{\sin{\frac{n}{2}\left[(\phi-Id)(\theta)-(\phi-Id)(\alpha)\right]}}
{\sin{\frac{\theta-\alpha}{2}}}\right)^2d\theta d\alpha.
\end{equation*}
\begin{equation*}
I_2^r = \cfrac{1}{2\pi} \iint_{|\theta-\alpha|\geq l} \left(
\cfrac{\sin{\frac{n}{2}\left[(\phi-Id)(\theta)-(\phi-Id)(\alpha)\right]}}
{\sin{\frac{\theta-\alpha}{2}}}\right)^2d\theta d\alpha.
\end{equation*}
Since for all $\delta' \in (0,2)$ and for small enough $l$, one has
\begin{equation*}
\left( \cfrac{\sin{\frac{n}{2}\left[(\phi-Id)(\theta)-(\phi-Id)(\alpha)\right]}}
{\sin{\frac{\theta-\alpha}{2}}}\right)^2 \leq C n^{\delta'}
\|\phi'-1\|_{\infty}^{\delta'-2}.
\end{equation*}
One checks that for $\delta' \in(1,2)$, there is a constant
$C(\delta')$ such that
\begin{equation*}
I_2^d \leq C(\delta') ~n^{\delta'}~
\|\phi'-1\|_{\infty}^{\delta'}.
\end{equation*}
Concerning $I_2^r$, we get easily for all $\delta''\in(0,1)$
\begin{equation*}
I_2^r \leq C~ n^{2\delta''}~ \|\phi'-1\|_{\infty}^{2\delta''}.
\end{equation*}
Summing up the estimates for $I_1$ and $I_2^d$ and $I_2^r$, we get
the stated result \eqref{inegalite:composition:atomique}.
\end{proofof}
Let us notice that the exponent $1+2\delta$ can hardly by reduced
since its main part $I_1$ is deduced from  the  Parseval's
equality. Another remark is that the constant $C(\delta)$ blows up
when $\delta \rightarrow 0$. We can now prove Theorem
\ref{theoreme:precomposition}.

\begin{proofof}{Theorem \ref{theoreme:precomposition}}
In a first step, we assume that $u$ is a trigonometric polynomial
namely
\begin{equation*}
u(\theta) = \sum_{|k|\leq n} c_k(u) e^{ik\theta}.
\end{equation*}
We fix $\delta'$ in $(0,\delta)$ and set $\alpha =
\delta/\delta'-1>0$. Then we have
\begin{eqnarray*}
\|u\circ\phi-u\|_{\sH^{1/2}(S^1)} &\leq &\sum_{|k|\leq n} |c_k(u)| ~\|
e^{ik\phi(\theta)}-e^{ik\theta} \|_{\sH^{1/2}(S^1)}\\
&\leq &C(\delta')~
\omega_{\delta'}(\|\phi-Id\|_{\sW^{1,\infty}(S^1)})~\sum_{|k|\leq n}
k^{1/2+\delta'}|c_k(u)|.
\end{eqnarray*}
Writing
\begin{equation*}
\cfrac{1}{2} +\delta' = -\cfrac{1}{2} -\alpha \delta' +1
+\delta'(1+\alpha) = -\cfrac{1}{2} -\alpha \delta' + 1+\delta,
\end{equation*}
we get
\begin{equation*}
\sum_{|k|\leq n} k^{1/2+\delta'}~|c_k(u)|~\leq \left( \sumZe |k|^{-1-2 \alpha \delta'}\right)^{1/2} ~
\left( \sumZe |k|^{2+2 (1+\alpha) \delta'}~|c_k(u)|^2 \right)^{1/2}.
\end{equation*}
Thanks to the Fatou's property that is satisfied by  the  Sobolev
spaces (\cite{Triebel}), we extend the result to the most general
Sobolev space  $\sH^{1+\delta}(S^1)$.
\end{proofof}

%%%%%%%%%%%%%%%%%%%%%%%%%%%%%%%%%%%%%%%%%%

\bibliographystyle{plain}

\end{document}